\documentclass[a4paper,12pt,twoside,leqno,final]{amsart}
\usepackage{amsmath}
\usepackage{amssymb}

\newtheorem{thm}{Theorem}[section]

\newtheorem{cor}[thm]{Corollary}

\newcommand{\C}{{\mathbb C}}
\newcommand{\D}{{\mathbb D}}
\newcommand{\R}{{\mathbb R}}
\newcommand{\T}{{\mathbb T}}
\newcommand{\Z}{{\mathbb Z}}
\newcommand{\N}{{\mathbb N}}
\newcommand{\cN}{{\mathcal N}}

\newcommand{\K}{{\mathcal{K}}}

\newcommand{\La}{\Lambda}

\newcommand{\f}{\frac}
\newcommand{\ov}{\overline}
\newcommand{\al}{\alpha}

\newcommand{\ga}{\gamma}

\newcommand{\de}{\delta}
\newcommand{\la}{\lambda}
\newcommand{\ze}{\zeta}
\renewcommand{\th}{\theta}
\newcommand{\si}{\sigma}
\newcommand{\ph}{\varphi}
\newcommand{\om}{\omega}
\newcommand{\Om}{\Omega}

\newcommand{\const}{\text{\rm const}}

\newcommand{\inn}{\text{\rm inn}}

\numberwithin{equation}{section}

\title[Boundary Gauss--Lucas type theorems on the disk]
{Boundary Gauss--Lucas type theorems \\ 
on the disk}
\author{Konstantin M. Dyakonov}
\address{BGSMath and Universitat de Barcelona, Departament de 
Matem\`atiques i Inform\`atica, Gran Via 585, E-08007 Barcelona, Spain}
\address{ICREA, Pg. Llu\'is Companys 23, E-08010 Barcelona, Spain}
\email{konstantin.dyakonov@icrea.cat}
\thanks{Supported in part by grant MTM2014-51834-P from El Ministerio de Econom\'ia 
y Competitividad (Spain) and grant 2014-SGR-289 from AGAUR (Generalitat de Catalunya).}

\begin{document}
\begin{abstract}
The classical Gauss--Lucas theorem describes the location of the critical points of a polynomial. There is also a hyperbolic version, due to Walsh, in which the role of polynomials is played by finite Blaschke products on the unit disk. We consider similar phenomena for generic inner functions, as well as for certain \lq\lq locally inner" self-maps 
of the disk. More precisely, we look at a unit-norm function $f\in H^\infty$ that has an angular derivative on a set 
of positive measure (on the boundary) and we assume that its inner factor, $I$, is nontrivial. Under certain 
conditions to be discussed, it follows that $f'$ must also have a nontrivial inner factor, say $J$, and we study 
the relationship between the boundary singularities of $I$ and $J$. Examples are furnished to show that our 
sufficient conditions cannot be substantially relaxed.
\end{abstract}

\maketitle

\section{Introduction and results} 

The functions considered in this paper are holomorphic self-maps of the unit disk. Our purpose is to find out 
when the presence of a nontrivial inner factor in the function's canonical factorization (i.e., the property 
of being non-outer) survives differentiation. This clearly happens if the original function, $f$, has multiple 
zeros, since these will also be zeros for $f'$. Somewhat less obvious is the fact that, under certain natural 
assumptions, the passage from $f$ to $f'$ preserves singular inner factors. (In a sense, these are responsible 
for the boundary zeros of infinite multiplicity.) Much subtler is the case of simple zeros inside, and this is 
what chiefly interests us here. The results that arise can be viewed as descendants of the classical Gauss--Lucas 
theorem, or rather of its disk version, and we begin by recalling those prototypical results. 
\par If $f$ is a holomorphic -- or perhaps meromorphic -- function living (at least) on a domain $\Om\subset\C$, 
we write $\mathcal Z_\Om(f)$ for its zero set there. Thus, 
$$\mathcal Z_\Om(f):=\{z\in\Om:\,f(z)=0\}.$$ 
The Gauss--Lucas theorem tells us that, given a nonconstant polynomial $P$, the set of its critical points, 
$\mathcal Z_\C(P')$, is contained in the convex hull of $\mathcal Z_\C(P)$; see, e.g., \cite[Chapter 2]{ShSm}. 
\par This fact admits a certain \lq\lq hyperbolic" analogue, which was found by Walsh \cite{W}. The plane $\C$ is 
now replaced by the disk 
$$\D:=\{z\in\C:|z|<1\},$$ 
while the role of polynomials is played by {\it finite Blaschke products}. A (finite) Blaschke product $B$ of 
degree $n$ is, by definition, given by 
$$B(z)=c\prod_{j=1}^n\f{z-a_j}{1-\bar a_jz}$$ 
(with some $a_1,\dots,a_n\in\D$ and a unimodular constant $c$), a formula known to provide the general form of 
an $n$-to-$1$ mapping from $\D$ onto itself. Now, Walsh's theorem says that, for such a $B$, the set 
$\mathcal Z_\D(B')$ is contained in the hyperbolic convex hull of $\mathcal Z_\D(B)$, defined appropriately; 
see \cite{W} for a precise statement. 
\par It should be mentioned that the set $\mathcal Z_\C(P')$ in the Gauss--Lucas theorem is 
automatically -- and trivially -- nonempty, provided that $\deg{P}\ge2$. Similarly, in Walsh's theorem, we have 
$\mathcal Z_\D(B')\ne\emptyset$ whenever $n$, the degree of $B$, satisfies $2\le n<\infty$. To see why, assume 
that $B'$ does not vanish on the set $\{0,a_1,\dots,a_n\}$ (otherwise the statement is trivial) and note that 
$\mathcal Z_\C(B')=\mathcal Z_\C(B'/B)$. The formula 
$$\f{B'(z)}{B(z)}=\sum_{j=1}^n\f{1-|a_j|^2}{(z-a_j)(1-\bar a_jz)}$$ 
shows then that the (nonempty) set $\mathcal Z_\C(B')$ is symmetric with respect to the 
circle $\T:=\partial\D$, so precisely one half of its points must be in $\D$. 
\par Thus, the Gauss--Lucas and Walsh theorems actually assert the {\it existence} of critical points in the 
appropriate region and also describe their {\it location}; this last part roughly amounts to saying that the 
zeros of $P'$ or $B'$ are to be found not too far from those of $P$ or $B$, respectively. 
\par Our purpose is to elaborate on Walsh's theorem by moving from finite Blaschke products to infinite 
ones, as well as to generic inner functions, and still further -- namely, to fairly general 
analytic self-maps of the disk -- and to study similar (Gauss--Lucas type) phenomena in these cases. 
\par At this point, we pause to recall some basic terminology and notation. A function $\th$ in $H^\infty$ 
(i.e., a bounded holomorphic function on $\D$) is said to be {\it inner} if $\lim_{r\to1^-}|\th(r\ze)|=1$ for 
$m$-almost all $\ze\in\T$. Here and throughout, $m$ is the normalized Lebesgue measure on the 
unit circle $\T$, so that $dm(\ze)=(2\pi)^{-1}|d\ze|$. 
It is well known that every inner function $\th$ can be factored canonically as 
$\th=\la BS$, where $\la\in\T$ is a constant, $B$ is a {\it Blaschke product} and $S$ is a {\it singular 
inner function}; see \cite[Chapter II]{G}. More explicitly, the factors involved are of the form 
\begin{equation}\label{eqn:infblaschke}
B(z)=B_{\{a_j\}}(z):=\prod_j\f{|a_j|}{a_j}\f{a_j-z}{1-\bar a_jz},\qquad z\in\D, 
\end{equation}
where $\{a_j\}\subset\D$ is a sequence -- possibly finite or empty -- with $\sum_j(1-|a_j|)<\infty$ 
(if $a_j=0$, one puts $|a_j|/a_j=-1$), and 
\begin{equation}\label{eqn:singinner}
S(z)=S_\mu(z):=\exp\left\{-\int_\T\f{\ze+z}{\ze-z}
d\mu(\ze)\right\},\qquad z\in\D,
\end{equation}
with $\mu$ a (nonnegative) singular measure on $\T$. The set 
$\T\cap\text{\rm clos}\left(\{a_j\}\cup\text{\rm supp}\,\mu\right)$ coincides with the {\it boundary spectrum} 
$\si(\th)$ of $\th$, defined as the set of its boundary singularities (i.e., the smallest closed set $E\subset\T$ 
such that $\th$ is analytic across $\T\setminus E$). 
\par Further, a zero-free holomorphic function $F$ on $\D$ is said to be {\it outer} if $\log|F|$ coincides 
with the harmonic extension (Poisson integral) of an integrable function on $\T$. When normalized by the 
condition $F(0)>0$, an outer function $F$ takes the form 
\begin{equation}\label{eqn:outer}
F(z)=\mathcal O_h(z):=\exp\left\{\int_\T\f{\ze+z}{\ze-z}\log h(\ze)\,dm(\ze)\right\},\qquad z\in\D,
\end{equation}
where $h$ is a nonnegative function on $\T$ with $\log h\in L^1(\T,m)$. This $h$ actually agrees with the 
nontangential boundary values of $|F|$ almost everywhere on $\T$. 
\par The functions $f$ that admit a factorization of the form $f=\th F$, with $\th$ inner and $F$ outer, are 
precisely those lying in the {\it Smirnov class} $\cN^+$; see \cite[Chapter II]{G}. Alternatively, we can define 
(or characterize) $\cN^+$ as the set of ratios $u/v$, where $u,v\in H^\infty$ and $v$ is outer. When $v$ is merely 
assumed to be zero-free on $\D$, such ratios range over the {\it Nevanlinna class} $\cN$. 
\par We write the canonical factorization of a function $f\in\cN^+$, $f\not\equiv0$, in the form 
\begin{equation}\label{eqn:canfact}
f=BSF,
\end{equation}
the three factors on the right being \eqref{eqn:infblaschke}, \eqref{eqn:singinner} and \eqref{eqn:outer}, 
respectively (and we take the liberty to ignore the unimodular constant factor involved). In particular, this 
canonical representation applies whenever $f$ is in the Hardy space $H^p$ with some $p\in(0,\infty]$ (see 
\cite[Chapter II]{G}); in fact, we have $H^p=\cN^+\cap L^p(\T,m)$. 

\par Now, going back to Walsh's theorem and trying to adapt it to an {\it infinite} 
Blaschke product $B$ (to begin with), we already have to face the new phenomenon that 
{\it the set $\mathcal Z_\D(B')$ may be empty}. Moreover, this may well happen for a Blaschke 
product $B$ with $B'\in\cN^+$. An example can be furnished as follows: fix 
a number $\al\in\D\setminus\{0\}$ and put 
$$B_\al(z):=\f{S(z)-\al}{1-\ov\al S(z)},$$ 
where $S$ is the \lq\lq atomic" singular inner function given by 
$$S(z):=\exp\left(\f{z+1}{z-1}\right).$$
It is well known (and easy to verify) that $B_\al$ is a Blaschke product. At the same time, 
differentiation yields 
$$B'_\al(z)=\f{1-|\al|^2}{\left(1-\bar\al S(z)\right)^2}\cdot\f{-2S(z)}{(z-1)^2},$$ 
and it is clear that the right-hand side is zero-free in $\D$. Furthermore, $B'_\al\in H^p$ for every 
$p\in(0,\f12)$, and the inner factor of $B'_\al$ is $S$. 

\par On the other hand, given a nonconstant inner function $\th$ with $\th'\in\cN^+$, it turns out that $\th'$ 
{\it must have a nontrivial inner factor, unless $\th$ is a M\"obius transformation} (see \cite[Corollary 2.2]{DCMFT} 
or \cite{DCR}). These observations seem to suggest that we modify our viewpoint appropriately. Namely, as long 
as our variations on Walsh's theme involving an inner function $\th$ are supposed to deal with something 
{\it a priori} existent, rather than pertain to the \lq\lq theory of the empty set", we feel that we 
should look at $\inn(\th')$, the inner factor of $\th'$, rather than at the zero set $\mathcal Z_\D(\th')$. 
(Here and below, we use the notation $\inn(f)$ for the inner factor of a function $f\in\cN^+$.) 
We should then try to understand the relationship between the (suitably defined) {\it smallness set} 
of $\inn(\th')$ and that of $\th$. More precisely, we shall be actually concerned with the boundary 
spectra $\si(\inn(\th'))$ and $\si(\th)$, i.e., with those parts of the unit circle $\T$ where the two smallness 
sets hit it. One consequence of our results is that 
$$\si(\inn(\th'))=\si(\th),$$ 
and we regard this as a boundary version of the Gauss--Lucas--Walsh theorem for inner functions. 
\par In fact, we are not going to restrict ourselves to inner functions, even though moving beyond this 
class makes things more complicated. This time, turning to a general function $f\in H^\infty$ with 
$f'\in\cN^+$, we can no longer expect that $f'$ will necessarily have an inner factor whenever $f$ does. For 
instance, suppose that $h$ is a holomorphic function on $\D$, with $\text{\rm Re}\,h$ bounded above, 
whose range $h(\D)$ contains infinitely many points of the form $c+2\pi ik$, where $c\in\C$ is fixed and 
$k$ ranges over (a subset of) $\Z$. The function $f:=e^h-e^c(\in H^\infty)$ then vanishes on the set 
$\bigcup_kh^{-1}(c+2\pi ik)$; thus, $f$ is divisible by an infinite Blaschke product, while $f'=h'e^h$ may 
well be outer. More sophisticated examples in this vein are given in Section 4 below. 
\par At the same time, we single out a class of $H^\infty$-functions that does obey the 
\lq\lq Gauss--Lucas principle", or perhaps the \lq\lq Walsh principle", meaning that the property of 
being non-outer is inherited by $f'$ from $f$ and that the boundary spectra of the two inner factors, 
$\inn(f')$ and $\inn(f)$, are related appropriately. The class in question appears to be (almost) optimal. 

\par Before stating the results, let us recall that a unit-norm function $f\in H^\infty$ is said to 
possess an {\it angular derivative} (in the sense of Carath\'eodory) at a point $\ze\in\T$ if 
both $f$ and $f'$ have nontangential limits at $\ze$ and, once we agree to denote the two limits 
by $f(\ze)$ and $f'(\ze)$, the former of these satisfies $|f(\ze)|=1$. The classical Julia--Carath\'eodory 
theorem (see \cite[Chapter VI]{B}, \cite[Chapter I]{Car} or \cite[Chapter VI]{Sar2}) 
asserts that this happens if and only if 
$$\liminf_{z\to\ze}\f{1-|f(z)|}{1-|z|}<\infty.$$
Further, given a point $z\in\D$, we shall denote by $\om_z$ the {\it harmonic measure} associated with it. 
Thus $d\om_z=P_zdm$ on $\T$, where $P_z$ stands for the corresponding {\it Poisson kernel}: 
$$P_z(\ze):=\f{1-|z|^2}{|\ze-z|^2},\qquad\ze\in\T.$$ 
The quantity $\om_z(E)=\int_Ed\om_z$, where $E$ is a (Lebesgue) measurable subset of $\T$, can be roughly 
interpreted as the normalized angle at which $E$ is seen from $z$. 
\par Also, we need to recall that a Blaschke product $b$ with zeros $\{z_j\}$ is said to be {\it thin} if 
$$\lim_{k\to\infty}\prod_{j:\,j\ne k}\left|\f{z_j-z_k}{1-\bar z_jz_k}\right|=1,$$
a condition that can be rewritten in the form 
$$\lim_{k\to\infty}|b'(z_k)|\left(1-|z_k|^2\right)=1.$$
The sequence $\{z_j\}$ itself is then also called thin, whereas non-thin sequences (and the corresponding 
Blaschke products) will be termed {\it thick}. In the literature, one encounters thin (or thick) sequences in 
many places. In particular, they turn up in connection with maximal ideals in uniform algebras and with various 
interpolation problems. One of the first occurrences can be found in \cite{Wo}; see also \cite{DN, GM, SW}. 

\par Now suppose that $\mathcal E$ is a (Lebesgue) measurable subset of $\T$, and 
$$\widetilde{\mathcal E}:=\T\setminus\mathcal E$$ 
(this notation will be used throughout), while $f$ is an $H^\infty$-function with $f'\in\cN^+$. 
Further, let $\inn(f)=BS$, where $B$ is a Blaschke product and $S$ a singular inner function. 
We then write 
$$\si^{\rm i}_{\mathcal E}(B):=\si(B)\cap\text{\rm ess\,int\,}\mathcal E,$$
where $\text{\rm ess\,int\,}\mathcal E$ is the {\it essential interior} of $\mathcal E$. (By definition, 
a point $\ze\in\T$ is in $\text{\rm ess\,int\,}\mathcal E$ if there exists a set $\La\subset\T$ with $m(\La)=0$ 
such that $\ze$ is an interior point of $\mathcal E\cup\La$ with respect to $\T$.) Also, we shall denote 
by $\si^{\rm b}_{\mathcal E,f}(B)$ the set of points $\ze\in\T\setminus\text{\rm ess\,int\,}\mathcal E$ 
with the following property: there exists a thick sequence $\{z_n\}\subset\mathcal Z_\D(B)$ with 
$z_n\to\ze$ satisfying 
\begin{equation}\label{eqn:firstcond}
\om_{z_n}(\widetilde{\mathcal E})\log\f1{1-|z_n|}\to0
\end{equation}
and 
\begin{equation}\label{eqn:secondcond}
\int_{\widetilde{\mathcal E}}\log|f'|\,d\om_{z_n}\to0. 
\end{equation}
The superscripts \lq\lq i" and \lq\lq b" in $\si^{\rm i}_{\mathcal E}(B)$ and $\si^{\rm b}_{\mathcal E,f}(B)$ 
stand for \lq\lq interior" and \lq\lq boundary", respectively. (It should be noted that the latter set is 
contained in $\text{\rm clos\,}\mathcal E$, so its elements are \lq\lq essentially boundary" points for 
$\mathcal E$.) Finally, we put 
$$\si_{\mathcal E}(f):=\si(S)\cup\si^{\rm i}_{\mathcal E}(B)\cup\si^{\rm b}_{\mathcal E,f}(B).$$ 

\begin{thm}\label{thm:newmainres} Let $f\in H^\infty$ be a nonconstant function with $\|f\|_\infty=1$, 
and let $\mathcal E$ be a measurable subset of $\T$ such that $f$ has an angular derivative almost everywhere 
on $\mathcal E$. Suppose that each of the three factors in the canonical factorization \eqref{eqn:canfact} has 
its derivative in $\cN^+$ (whence also $f'\in\cN^+$). Assume, finally, that $\si_{\mathcal E}(f)\ne\emptyset$. 
Then $f'$ has a nontrivial inner factor, say $J$, with $\si(J)\supset\si_{\mathcal E}(f)$. 
\end{thm}

\par The last inclusion should be compared with the fact that $\si(J)$ is always contained in $\si(f)$, 
the set of boundary singularities for $f$. Indeed, if $f$ is analytic in a neighborhood of a point $\ze_0\in\T$, 
then so is $f'$, and hence also its inner factor, $J$ (see \cite[Chapter II]{G} in connection with the latter 
implication). 

\par Also, in the theorem above, we may replace the hypothesis that $f$ has an angular derivative a.\,e. on 
$\mathcal E$ by the seemingly weaker condition that $|f|=1$ a.\,e. on $\mathcal E$. (The reason is that the other 
assumptions imply the existence of nontangential limits for $f'$ a.\,e. on $\T$.) The functions $f$ that arise 
can thus be viewed as \lq\lq locally inner". Because of the role that angular derivatives play in Theorem 
\ref{thm:arcad} below, we have chosen to state Theorem \ref{thm:newmainres} in similar terms; the relation 
between the two results might in this way become clearer. 

\par We now make a remark concerning the meaning of conditions \eqref{eqn:firstcond} and \eqref{eqn:secondcond} 
that were used to define the set $\si^{\rm b}_{\mathcal E,f}(B)$. Given an (essentially) boundary 
point $\ze$ of $\mathcal E$ and a sequence $\{z_n\}\subset\mathcal Z_\D(B)$ with $z_n\to\ze$, the two 
conditions basically mean that the $z_n$'s tend to $\ze$ tangentially enough \lq\lq on the $\mathcal E$ side" 
(i.e., they lie much closer to $\mathcal E$ than to $\widetilde{\mathcal E}$). The examples constructed at 
the end of the paper will show that nontangential convergence would not do, and moreover, that the 
qualitative tangency conditions \eqref{eqn:firstcond} and \eqref{eqn:secondcond} cannot be substantially 
relaxed. 

\par One may find it unfortunate that condition \eqref{eqn:secondcond} involves $f'$, instead of being stated in 
terms of $f$ alone. We note, however, that it only depends on the boundary values of $|f'|$ (or, equivalently, on 
the {\it outer} factor of $f'$), whereas the conclusion of Theorem \ref{thm:newmainres} concerns the {\it inner} 
factor of $f'$. Besides, under further hypotheses, we shall come up with simpler sufficient conditions replacing 
\eqref{eqn:firstcond} and \eqref{eqn:secondcond} that will lead to more transparent formulations. 

\par The following corollary deals with the situation where $\si(B)$ is contained 
in $\text{\rm ess\,int\,}\mathcal E$, in which case we have $\si^{\rm i}_{\mathcal E}(B)=\si(B)$, 
$\si^{\rm b}_{\mathcal E,f}(B)=\emptyset$ and 
$$\si_{\mathcal E}(f):=\si(B)\cup\si(S)=\si(BS).$$  

\begin{cor}\label{cor:innoutthick} Let $f\in H^\infty$ be a nonconstant function with $\|f\|_\infty=1$, 
and let $\mathcal E$ be a measurable subset of $\T$ such that $f$ has an angular derivative almost everywhere 
on $\mathcal E$. Suppose that each of the three factors in the canonical factorization \eqref{eqn:canfact} has 
its derivative in $\cN^+$. Assume, finally, that $\si(BS)\ne\emptyset$, while 
$\si(B)\subset\text{\rm ess\,int\,}\mathcal E$. Then $J:=\inn(f')$ is a nontrivial inner function 
and $\si(BS)\subset\si(J)$. 
\end{cor}

\par In the special case where $\mathcal E=\T$, this reduces to the following result. 

\begin{cor}\label{cor:innthick} 
Let $\th$ be a nonconstant inner function, other than a M\"obius transformation, with $\th'\in\cN^+$. Then 
$\mathcal J:=\inn(\th')$ is a nontrivial inner function and $\si(\th)=\si(\mathcal J)$. 
\end{cor} 

\par In view of the discussion following Theorem \ref{thm:newmainres}, we have the (trivial) inclusion 
$\si(\mathcal J)\subset\si(\th)$. Now, if $\si(\th)\ne\emptyset$, Corollary \ref{cor:innthick} is a special 
case of the preceding result (just take $\mathcal E=\T$ and $f=\th=BS$). Otherwise, we are only concerned with 
the nontriviality of $\inn(\th')$, and this is guaranteed by the above-mentioned result from \cite{DCMFT, DCR}. 

\par Regarding the hypothesis $\th'\in\cN^+$ (for $\th$ inner), we recall that this is actually equivalent 
to $\th'\in\cN$. Furthermore, each of these holds if and only if $\log^+|\th'|\in L^1(\T,m)$, where $\th'$ 
is understood as the angular derivative. Also, for $\th=BS$ to satisfy $\th'\in\cN$ (or $\th'\in\cN^+$), 
it is necessary and sufficient that both $B'$ and $S'$ be in $\cN$ (or $\cN^+$). These results are due to Ahern 
and Clark; see \cite[Corollary 4]{AC}. 

\par The next result, also a consequence of Theorem \ref{thm:newmainres}, contains a simple sufficient 
condition for $\inn(f')$ to be nontrivial when $\mathcal E$ is taken to be an arc. 
In what follows, we write $\mathcal A$ for the {\it disk algebra} $H^\infty\cap C(\T)$. 

\begin{thm}\label{thm:arcad} Let $\mathcal E=\{e^{it}:\,0\le t\le t_0\}$, where $0<t_0\le\pi$, and let 
$F\in H^\infty$ be an outer function such that $\|F\|_\infty=1$, $F'\in\mathcal A$, and $|F|=1$ on $\mathcal E$. 
Further, suppose $\{z_n\}\subset\D$ is a thick sequence with the properties that $\text{\rm Im}\,z_n>0$ ($n\in\N$), 
$\lim_{n\to\infty}z_n=1$ and 
\begin{equation}\label{eqn:angderbla}
\sum_n\f{1-|z_n|^2}{|1-z_n|^2}<\infty.
\end{equation}
Finally, assume that the Blaschke product $B=B_{\{z_n\}}$ satisfies $B'\in\cN^+$ and put $f:=BF$. Then $f'$ 
lies in $\cN^+$ and has a nontrivial inner factor, $J$, with $1\in\si(J)$. 
\end{thm}

\par A few remarks are in order. First, it is easy to construct an outer function $F$ satisfying the hypotheses 
of Theorem \ref{thm:arcad} by defining its modulus $|F|\big|_\T=:h$ appropriately. Namely, it suffices to assume 
that $h\in C^{2+\varepsilon}(\T)$ for some $\varepsilon>0$ (i.e., that $h''$ is Lipschitz continuous of order 
$\varepsilon$), in addition to the obvious conditions that $\log h\in L^1(\T,m)$, $0\le h\le1$ on $\T$, 
and $h\big|_\mathcal E=1$. Now, for $F=\mathcal O_h$, the fact that $F'\in\mathcal A$ (and actually the stronger 
conclusion that $F\in C^{1+\varepsilon/2}(\T)$) is guaranteed by the Havin--Shamoyan(--Carleson--Jacobs) 
theorem; see \cite{DActa, H, Shi}. In the case where $h$ is strictly positive, the regularity assumption 
can be relaxed to $h\in C^{1+\varepsilon}(\T)$, since $F=\mathcal O_h$ will then be in the same class; this 
follows from standard properties of the Hilbert transform, see \cite[Chapter III]{G}. 

\par Secondly, condition \eqref{eqn:angderbla} means precisely that $B$ has an angular derivative at the point 
$1$. Roughly speaking, it says that the zero sequence $\{z_n\}$ approaches its limit point $1$ in a suitably 
tangential manner. The (sufficient) tangency condition expressed by \eqref{eqn:angderbla} should be compared 
with the weaker condition \eqref{eqn:firstcond}, which, alone, does not suffice to conclude that 
$\inn(f')$ is nontrivial; see Example 2 in Section 4 below. 

\par Thirdly, we could have stated Theorem \ref{thm:arcad} in a more general form, where a nontrivial singular 
factor $S$ (with $S'\in\cN^+$) is present in the canonical factorization \eqref{eqn:canfact}. The conclusion 
would have been that $\si(S)\cup\{1\}\subset\si(J)$. However, the main issue being the location of zeros, 
we have chosen to restrict ourselves to the current version. 

\par The proofs of Theorems \ref{thm:newmainres} and \ref{thm:arcad} are given in Sections 2 and 3, 
respectively, while the last section contains a couple of examples to the effect that the hypotheses of 
Theorem \ref{thm:newmainres} are close to being sharp. 

\par I thank Pascal Thomas for a helpful remark concerning the formulation of Theorem \ref{thm:newmainres}. 

\medskip

\section{Proof of Theorem \ref{thm:newmainres}} 

First of all, since 
\begin{equation}\label{eqn:diffthree}
f'=B'SF+BS'F+BSF',
\end{equation}
our hypotheses on the three factors guarantee that $f'\in\cN^+$. 
\par Furthermore, because $\si_{\mathcal E}(f)\ne\emptyset$, we know that either $\si(S)\ne\emptyset$ 
or 
$$\si^{\rm i}_{\mathcal E}(B)\cup\si^{\rm b}_{\mathcal E,f}(B)=:\si_{\mathcal E,f}(B)\ne\emptyset.$$ 
Assuming that $\si(S)\ne\emptyset$ (so that the singular factor 
$S$ is nontrivial), we now rewrite \eqref{eqn:diffthree} as 
$$\f{f'}S=B'F+BF'+BF\f{S'}S$$ 
and claim that each of the three terms on the right is in $\cN^+$. Indeed, for the last term, this is 
ensured by Ahern and Clark's result (see \cite[Corollary 4]{AC}) which says that $S'/S\in\cN^+$ whenever 
$S'\in\cN^+$; the preceding terms present no difficulty. It follows that $f'/S\in\cN^+$, and so $J:=\inn(f')$ 
is divisible by $S$. In particular, we have then $J\not\equiv\const$ and $\si(S)\subset\si(J)$. 
\par To deal with the case where $\si_{\mathcal E,f}(B)\ne\emptyset$, more work is needed. Let $G$ stand 
for the outer factor of $f'$, so that 
$$G(z)=\exp\left\{\int_\T\f{\ze+z}{\ze-z}\log |f'(\ze)|\,dm(\ze)\right\},\qquad z\in\D,$$ 
and let $G_{\mathcal E}$ (resp., $G_{\widetilde{\mathcal E}}$) be defined by a similar formula, where the 
integral is taken over $\mathcal E$ (resp., over $\widetilde{\mathcal E}:=\T\setminus\mathcal E$). Thus, 
in particular, $G_{\mathcal E}$ is the outer function with 
modulus $|f'|\chi_{\mathcal E}+\chi_{\widetilde{\mathcal E}}$ and 
$G_{\widetilde{\mathcal E}}=G/G_{\mathcal E}$. 
\par Our plan is to deduce the nontriviality of the inner function $J:=f'/G$, plus the fact that 
$\si(J)$ contains $\si_{\mathcal E,f}(B)$, from the inequality 
\begin{equation}\label{eqn:crucineq}
|J(z)|\left|G_{\widetilde{\mathcal E}}(z)\right|\le\f{|f'(z)|\,(1-|z|^2)}{1-|f(z)|^2}\cdot\ga_\mathcal E(z),
\qquad z\in\D,
\end{equation}
where 
$$\ga_\mathcal E(z):=\left\{\f2{(1-|z|)\,\om_z(\widetilde{\mathcal E})}\right\}^{\om_z(\widetilde{\mathcal E})}.$$
This crucial estimate will be established later on. Right now, we take it for granted and complete the proof. 

\par Suppose $\ze\in\si^{\rm i}_{\mathcal E}(B)$. By adding a suitable null-set to $\mathcal E$ if necessary, we 
may assume that $\ze$ is an interior point of $\mathcal E$. To show that $\ze\in\si(J)$, we argue by contradiction. 
Let $\Gamma$ be an open subarc of $\T$ with $\ze\in\Gamma\subset\mathcal E$ such that $J$ is analytic 
across $\Gamma$, and fix a point $\xi\in\Gamma$. Then 
\begin{equation}\label{eqn:jtoone}
|J(z)|\to1\quad\text{\rm as}\quad z\to\xi
\end{equation}
(it is always understood that $z$ is restricted to $\D$). Furthermore, since $\xi$ lies at a positive distance 
from $\widetilde{\mathcal E}$, the Poisson kernels $P_z$ satisfy 
$$\sup\{P_z(\eta):\,\eta\in\widetilde{\mathcal E}\}\le C\cdot(1-|z|)$$ 
whenever $z$ is close enough to $\xi$; here $C=C(\xi,\mathcal E)$ is a positive constant. This implies 
that the quantities 
$$\log\ga_\mathcal E(z)=\om_{z}(\widetilde{\mathcal E})\cdot
\log\left\{\f2{(1-|z|)\,\om_{z}(\widetilde{\mathcal E})}\right\}$$ 
and 
$$\log\left|G_{\widetilde{\mathcal E}}(z)\right|=\int_{\widetilde{\mathcal E}}\log|f'(\eta)|\,d\om_{z}(\eta)$$ 
both tend to $0$ as $z\to\xi$, whence 
\begin{equation}\label{eqn:pugalo}
\ga_\mathcal E(z)\to1\quad\text{\rm and}\quad
\left|G_{\widetilde{\mathcal E}}(z)\right|\to1,\quad\text{\rm as }z\to\xi. 
\end{equation}
Combining \eqref{eqn:crucineq} with \eqref{eqn:jtoone} and \eqref{eqn:pugalo}, in conjunction with the 
Schwarz--Pick estimate 
$$\f{|f'(z)|\,(1-|z|^2)}{1-|f(z)|^2}\le1,$$
we see that 
\begin{equation}\label{eqn:solonka}
\f{|f'(z)|\,(1-|z|^2)}{1-|f(z)|^2}\to1\quad\text{\rm as}\quad z\to\xi.
\end{equation}
This being true for every $\xi\in\Gamma$, we invoke a result of Kraus, Roth and Ruscheweyh 
(see \cite[Theorem 1.1]{KRR}) to conclude from \eqref{eqn:solonka} that $f$ is analytic across $\Gamma$. 
However, this conclusion is incompatible with the fact that the zeros of $f$ cluster at $\ze(\in\Gamma)$. 
The contradiction proves that $\ze\in\si(J)$. 
\par We have thus established the inclusion $\si^{\rm i}_{\mathcal E}(B)\subset\si(J)$. In particular, 
it follows that $J$ is nonconstant whenever $\si^{\rm i}_{\mathcal E}(B)\ne\emptyset$. 

\par Now suppose that $\ze\in\si^{\rm b}_{\mathcal E,f}(B)$. By definition, this means that we can find 
a thick sequence $\{z_n\}\subset\mathcal Z_\D(B)$ with $z_n\to\ze$ satisfying \eqref{eqn:firstcond} 
and \eqref{eqn:secondcond}. Let $b$ be the Blaschke product with zeros $\{z_n\}$, so that $b$ is 
a subproduct of $B$ and 
\begin{equation}\label{eqn:liminfbla}
\liminf_{n\to\infty}|b'(z_n)|\left(1-|z_n|^2\right)<1.
\end{equation}
Also, we have $f=gb$ for some $g\in H^\infty$ with $\|g\|_\infty=1$; this in turn implies that 
$$|f'(z_n)|=|g(z_n)|\cdot|b'(z_n)|\le|b'(z_n)|$$ 
for each $n$. Consequently, applying \eqref{eqn:crucineq} with $z=z_n$ yields 
\begin{equation}\label{eqn:applcru}
|J(z_n)|\left|G_{\widetilde{\mathcal E}}(z_n)\right|
\le\ga_\mathcal E(z_n)\cdot|b'(z_n)|\cdot\left(1-|z_n|^2\right).
\end{equation}
Now, \eqref{eqn:firstcond} shows that the quantity 
$$\log\ga_\mathcal E(z_n)=\om_{z_n}(\widetilde{\mathcal E})\cdot
\log\left\{\f2{(1-|z_n|)\,\om_{z_n}(\widetilde{\mathcal E})}\right\}$$ 
tends to $0$ as $n\to\infty$, while \eqref{eqn:secondcond} leads us to a similar conclusion about the quantity 
$$\log\left|G_{\widetilde{\mathcal E}}(z_n)\right|=\int_{\widetilde{\mathcal E}}\log|f'(\eta)|\,d\om_{z_n}(\eta).$$ 
Therefore, 
\begin{equation}\label{eqn:limgag}
\lim_{n\to\infty}\ga_\mathcal E(z_n)=1\qquad\text{\rm and}\qquad
\lim_{n\to\infty}\left|G_{\widetilde{\mathcal E}}(z_n)\right|=1.
\end{equation}
Finally, taking \eqref{eqn:liminfbla} and \eqref{eqn:limgag} into account, we deduce from \eqref{eqn:applcru} that 
$$\liminf_{n\to\infty}|J(z_n)|<1.$$ 
Recalling that $z_n\to\ze(\in\T)$, we readily conclude that $\ze\in\si(J)$. 
\par Now we know that $\si^{\rm b}_{\mathcal E,f}(B)\subset\si(J)$; and this clearly implies that $J\not\equiv\const$ 
whenever $\si^{\rm b}_{\mathcal E,f}(B)\ne\emptyset$. 

\par It remains to verify \eqref{eqn:crucineq}. Let $z\in\D$ be fixed. Then, for almost all $\ze\in\mathcal E$, 
Julia's lemma (see \cite[p.\,41]{G}) gives 
\begin{equation}\label{eqn:julia}
\f{|f(\ze)-f(z)|^2}{1-|f(z)|^2}\le|f'(\ze)|\cdot\f{|\ze-z|^2}{1-|z|^2}, 
\end{equation}
or equivalently, 
\begin{equation}\label{eqn:juliabis}
\f{1-|z|^2}{1-|f(z)|^2}\cdot\left|\f{1-\ov{f(z)}f(\ze)}{1-\ov z\ze}\right|^2\le|f'(\ze)| 
\end{equation}
(recall that $|f(\ze)|=1$ whenever $f$ has an angular derivative at $\ze$). 
Next, we consider the $H^\infty$-function 
\begin{equation}\label{eqn:fzphi}
\Phi_z(w):=\f{1-|z|^2}{1-|f(z)|^2}\cdot\left(\f{1-\ov{f(z)}f(w)}{1-\ov zw}\right)^2
\end{equation}
and rewrite \eqref{eqn:juliabis} in the form 
\begin{equation}\label{eqn:estone}
|\Phi_z(\ze)|\le|f'(\ze)|,\qquad\ze\in\mathcal E.
\end{equation}
Further, we define $\Psi_z$ to be the outer function with modulus 
$$\left|\Psi_z(\ze)\right|=|f'(\ze)|\cdot\chi_\mathcal E(\ze)
+|\Phi_z(\ze)|\cdot\chi_{\widetilde{\mathcal E}}(\ze),\qquad \ze\in\T,$$ 
and observe that 
\begin{equation}\label{eqn:estoncirc}
|\Phi_z(\ze)|\le|\Psi_z(\ze)|,\qquad\ze\in\T.
\end{equation} 
In fact, for $\ze\in\mathcal E$, this inequality coincides with \eqref{eqn:estone}, while 
for $\ze\in\widetilde{\mathcal E}$ the two sides are obviously equal. 
\par Since $\Psi_z$ is outer, the estimate \eqref{eqn:estoncirc} extends into $\D$, so that 
$$|\Phi_z(w)|\le|\Psi_z(w)|,\qquad w\in\D.$$
In particular, this holds for $w=z$, whence 
\begin{equation}\label{eqn:estatz}
|\Phi_z(z)|\le|\Psi_z(z)|.
\end{equation}
A glance at \eqref{eqn:fzphi} reveals that 
\begin{equation}\label{eqn:lefthand}
|\Phi_z(z)|=\Phi_z(z)=\f{1-|f(z)|^2}{1-|z|^2},
\end{equation}
and we take further steps to estimate $|\Psi_z(z)|$. 
\par We have 
\begin{equation}\label{eqn:righthand}
\log|\Psi_z(z)|=\int_\T\log|\Psi_z(\ze)|\,d\om_z(\ze)=I_1(z)+I_2(z),
\end{equation}
where 
\begin{equation}\label{eqn:int1}
I_1(z):=\int_\mathcal E\log|f'(\ze)|\,d\om_z(\ze)=\log|G_{\mathcal E}(z)|
\end{equation}
and 
\begin{equation}\label{eqn:int2}
I_2(z):=\int_{\widetilde{\mathcal E}}\log|\Phi_z(\ze)|\,d\om_z(\ze).
\end{equation}
The arithmetic/geometric mean inequality yields 
\begin{equation}\label{eqn:conc}
\begin{aligned}
I_2(z)&=\om_z(\widetilde{\mathcal E})\cdot\int_{\widetilde{\mathcal E}}\log|\Phi_z(\ze)|\,
\f{d\om_z(\ze)}{\om_z(\widetilde{\mathcal E})}\\
&\le\om_z(\widetilde{\mathcal E})\cdot\log\left\{\f1{\om_z(\widetilde{\mathcal E})}
\int_{\widetilde{\mathcal E}}|\Phi_z(\ze)|\,d\om_z(\ze)\right\}\\
&\le\om_z(\widetilde{\mathcal E})\cdot\log\left\{\f1{\om_z(\widetilde{\mathcal E})}
\int_\T|\Phi_z(\ze)|\,d\om_z(\ze)\right\}.
\end{aligned}
\end{equation}
We proceed by noticing that, for almost all $\ze\in\T$, 
\begin{equation}\label{eqn:fzze}
\begin{aligned}
|\Phi_z(\ze)|&\le\f2{1-|z|}\cdot\f{|1-\ov{f(z)}f(\ze)|^2}{1-|f(z)|^2}\\ 
&=c_{f,z}\left\{1-2\,\text{\rm Re}\left(\ov{f(z)}f(\ze)\right)+|f(z)|^2|f(\ze)|^2\right\}, 
\end{aligned}
\end{equation}
where 
$$c_{f,z}:=\f2{(1-|z|)\left(1-|f(z)|^2\right)}.$$ 
Also, we introduce the (harmonic) function 
$$u_z(\ze):=1-2\,\text{\rm Re}\left(\ov{f(z)}f(\ze)\right)+|f(z)|^2$$ 
and go on to observe that 
$$|\Phi_z(\ze)|\le c_{f,z}\cdot u_z(\ze),\qquad\ze\in\T.$$ 
(This follows from \eqref{eqn:fzze} and the fact that $|f(\ze)|\le1$ on $\T$.) Consequently, 
\begin{equation}\label{eqn:taburetka}
\int_\T|\Phi_z(\ze)|\,d\om_z(\ze)\le c_{f,z}\int_\T u_z(\ze)\,d\om_z(\ze)=c_{f,z}\cdot u_z(z)=\f2{1-|z|}. 
\end{equation}
Plugging the resulting inequality from \eqref{eqn:taburetka} into \eqref{eqn:conc}, we now get 
\begin{equation}\label{eqn:int2fin}
I_2(z)\le\om_z(\widetilde{\mathcal E})\log\f2{(1-|z|)\cdot\om_z(\widetilde{\mathcal E})}.
\end{equation}
\par This done, we combine \eqref{eqn:righthand} with \eqref{eqn:int1} and \eqref{eqn:int2fin} to infer that 
$$\log|\Psi_z(z)|\le\log|G_{\mathcal E}(z)|
+\om_z(\widetilde{\mathcal E})\log\f2{(1-|z|)\cdot\om_z(\widetilde{\mathcal E})}$$
and hence 
\begin{equation}\label{eqn:unitaz}
|\Psi_z(z)|\le|G_{\mathcal E}(z)|\left\{\f2{(1-|z|)\cdot\om_z(\widetilde{\mathcal E})}
\right\}^{\om_z(\widetilde{\mathcal E})}=|G_{\mathcal E}(z)|\cdot\ga_\mathcal E(z).
\end{equation}
Since
$$G_{\mathcal E}(z)=\f{G(z)}{G_{\widetilde{\mathcal E}}(z)}=\f{f'(z)}{G_{\widetilde{\mathcal E}}(z)J(z)},$$
we may further rewrite \eqref{eqn:unitaz} as 
\begin{equation}\label{eqn:unitazbis}
|\Psi_z(z)|\le\left|\f{f'(z)}{G_{\widetilde{\mathcal E}}(z)J(z)}\right|\cdot\ga_\mathcal E(z).
\end{equation}
On the other hand, recalling \eqref{eqn:estatz} and \eqref{eqn:lefthand}, we see that 
\begin{equation}\label{eqn:unitaztri}
|\Psi_z(z)|\ge\f{1-|f(z)|^2}{1-|z|^2}.
\end{equation}
Finally, a juxtaposition of \eqref{eqn:unitazbis} and \eqref{eqn:unitaztri} yields 
$$\f{1-|f(z)|^2}{1-|z|^2}\le\left|\f{f'(z)}{G_{\widetilde{\mathcal E}}(z)J(z)}\right|\cdot\ga_\mathcal E(z),$$ 
which is precisely \eqref{eqn:crucineq}. The proof is therefore complete. 

\medskip

\section{Proof of Theorem \ref{thm:arcad}}

It suffices to check that the hypotheses of the current theorem imply those of Theorem \ref{thm:newmainres}, 
with $S\equiv1$, and that $\si_{\mathcal E}(f)=\si^{\rm b}_{\mathcal E,f}(B)=\{1\}$. An application 
of Theorem \ref{thm:newmainres} will then do the job. 
\par First of all, the assumptions on $B$ and $F$ guarantee that $f=BF$ possesses an angular derivative 
everywhere on $\mathcal E=\{e^{it}:\,0\le t\le t_0\}$. Indeed, each of the two factors enjoys a similar 
property there; in particular, \eqref{eqn:angderbla} tells us that $B$ has an angular derivative at 
the endpoint $1$. This said, it remains to verify conditions \eqref{eqn:firstcond} and \eqref{eqn:secondcond}, 
where $\widetilde{\mathcal E}$ is the arc complementary to $\mathcal E$. 
\par The verification of \eqref{eqn:firstcond} is straightforward. In fact, from \eqref{eqn:angderbla} it clearly 
follows that 
$$\f{(1-|z_n|)^{1/2}}{|1-z_n|}\to0.$$
It is also obvious that 
$$(1-|z_n|)^{1/2}\log\f1{1-|z_n|}\to0.$$
Consequently, the product of the two quantities, which is 
$$\f{1-|z_n|}{|1-z_n|}\log\f1{1-|z_n|},$$
also tends to $0$ as $n\to\infty$. Together with the elementary estimate 
$$\om_{z_n}\left(\widetilde{\mathcal E}\right)\asymp\f{1-|z_n|}{|1-z_n|},$$ 
this yields \eqref{eqn:firstcond}. 
\par Here and below, the sign $\asymp$ is used to mean that the quantities involved are comparable (i.e., their 
ratio lies between two positive constants). 
\par To check \eqref{eqn:secondcond}, we are going to prove the following claim: {\it There exist numbers 
$C>1$ and $\de\in(0,\f\pi2)$ such that} 
\begin{equation}\label{eqn:zhaba}
C^{-1}\le|f'(e^{it})|\le C\quad\text{\it whenever}\quad-\de<t<0.
\end{equation}
Once this is established, \eqref{eqn:secondcond} comes out easily. Indeed, on the arc 
\begin{equation}\label{eqn:gade}
\ga_\de:=\left\{e^{it}:\,-\de<t<0\right\}
\end{equation}
we have $-M\le\log|f'|\le M$ with $M:=\log C$, whence 
\begin{equation}\label{eqn:firstest}
\left|\int_{\ga_\de}\log|f'|\,d\om_{z_n}\right|\le M\om_{z_n}\left(\ga_\de\right)
\le M\om_{z_n}\left(\widetilde{\mathcal E}\right)\to0. 
\end{equation}
Now, for $\ze\in\widetilde{\mathcal E}\setminus\ga_\de$, we have 
$$P_{z_n}(\ze)\le\const\cdot(1-|z_n|)$$ 
(because the limit point $1$ of the $z_n$'s lies at a positive distance from 
$\widetilde{\mathcal E}\setminus\ga_\de$), and so 
\begin{equation}\label{eqn:secondest}
\left|\int_{\widetilde{\mathcal E}\setminus\ga_\de}\log|f'|\,d\om_{z_n}\right|
\le\const\cdot(1-|z_n|)\int_\T\left|\log|f'|\right|dm\to0.
\end{equation}
Combining \eqref{eqn:firstest} and \eqref{eqn:secondest}, we arrive at \eqref{eqn:secondcond}. 
\par We now turn to proving the claim above; see \eqref{eqn:zhaba} and the italicized text preceding it. 
Since $f'=B'F+BF'$, the right-hand inequality in \eqref{eqn:zhaba} will be established as soon as we 
show that $|B'(e^{it})|$ is bounded for $-\f\pi2\le t<0$. To this end, we write $z_n=r_ne^{i\ph_n}$ 
(with $r_n>0$ and $0<\ph_n<\pi$) and estimate the quantity 
$$|B'(\ze)|=\sum_n\f{1-r_n^2}{|\ze-z_n|^2}$$
at a point $\ze=e^{it}$ with $-\f\pi2\le t<0$. There is no loss of generality in assuming that 
$r_n\ge\f12$ and $0<\ph_n\le\f\pi2$, since this is true for all but finitely many $z_n$'s. We now combine 
the elementary inequalities 
\begin{equation}\label{eqn:elemineq}
\f2{\pi^2}\left[(1-r_n)^2+(\ph_n-t)^2\right]\le|\ze-z_n|^2\le(1-r_n)^2+(\ph_n-t)^2, 
\end{equation}
valid in this case, with the fact that 
\begin{equation}\label{eqn:elemineqbis}
|\ph_n-t|=\ph_n+|t|\ge\ph_n
\end{equation}
to infer that 
\begin{equation}\label{eqn:abracadabra}
\begin{aligned}
|B'(\ze)|&\le\f{\pi^2}2\sum_n\f{1-r_n^2}{(1-r_n)^2+(\ph_n-t)^2}\\
&\le\f{\pi^2}2\sum_n\f{1-r_n^2}{(1-r_n)^2+\ph_n^2}\\
&\le\f{\pi^2}2\sum_n\f{1-r_n^2}{|1-z_n|^2}.
\end{aligned}
\end{equation}
The last quantity being finite by \eqref{eqn:angderbla}, it follows that $|B'(e^{it})|$ is bounded 
for $-\f\pi2\le t<0$. Consequently, 
\begin{equation}\label{eqn:fprimebdd}
\sup\left\{|f'(e^{it})|:\,-\f\pi2\le t<0\right\}<\infty,
\end{equation}
which proves \lq\lq half" of \eqref{eqn:zhaba}. 
\par Moving on to the left-hand inequality in \eqref{eqn:zhaba}, we first note that the modulus of the 
(angular) derivative $F'(\ze)$ at a point $\ze\in\mathcal E$ coincides with the nontangential limit of 
$$Q_F(z):=\f{1-|F(z)|^2}{1-|z|^2}$$ 
as $z\to\ze$; this forms part of the Julia--Carath\'eodory theorem. Secondly, we recall that 
$$Q_F(z)\ge\f{1-|F(0)|}{1+|F(0)|}=:\eta(=\eta_F)>0,\qquad z\in\D,$$
a well-known consequence of Schwarz's lemma (see, e.g., \cite{Sar1}). Therefore, 
\begin{equation}\label{eqn:feta}
|F'(\ze)|\ge\eta,\qquad\ze\in\mathcal E.
\end{equation}
Now let $N\in\N$ be a number such that 
\begin{equation}\label{eqn:choiceofn}
\sum_{n=N+1}^\infty\f{1-|z_n|^2}{|1-z_n|^2}<\f\eta{2\pi^{2}},
\end{equation}
and let $B_0$ and $B_1$ be the Blaschke products with zero sets $\{z_n:\,1\le n\le N\}$ and $\{z_n:\,n>N\}$, 
respectively. Then put $G:=FB_0$, so that $f=GB_1$. 
\par Since $F$ and $B_0$ both have an angular derivative on $\mathcal E$, the same is true for $G$. Moreover, 
it follows (see \cite[Corollary 1]{AC}) that 
$$|G'(\ze)|=|F'(\ze)|+|B'_0(\ze)|,\qquad\ze\in\mathcal E.$$ 
Recalling \eqref{eqn:feta}, we see that $|G'|\ge|F'|\ge\eta$ on $\mathcal E$; and since $G'(=F'B_0+FB'_0)$ 
is continuous on $\T$, we can find a number $\de\in(0,\f\pi 2)$ such that 
\begin{equation}\label{eqn:gongade}
|G'(\ze)|\ge\f\eta2,\qquad\ze\in\ga_\de
\end{equation}
(here $\ga_\de$ is the arc defined by \eqref{eqn:gade}). Furthermore, because $|G|\le1=|B_1|$ on $\T$, 
we have 
\begin{equation}\label{eqn:fprimebelow}
|f'|\ge|G'B_1|-|GB'_1|\ge|G'|-|B'_1|
\end{equation}
there; in particular, this holds on $\ga_\de$. 
\par Finally, we estimate the quantity 
\begin{equation}\label{eqn:bone}
|B'_1(\ze)|=\sum_{n=N+1}^\infty\f{1-|z_n|^2}{|\ze-z_n|^2}
\end{equation}
at a point $\ze=e^{it}$ with $-\f\pi2\le t<0$. As before, we write $z_n=r_ne^{i\ph_n}$, assuming that 
$r_n\ge\f12$ and $0<\ph_n\le\f\pi2$ (this is certainly true for $n>N$, with $N$ large enough), and we 
employ the elementary inequalities \eqref{eqn:elemineq} and \eqref{eqn:elemineqbis} to estimate the sum 
in \eqref{eqn:bone}. The estimate, which mimics \eqref{eqn:abracadabra}, reads 
\begin{equation}\label{eqn:abra}
\begin{aligned}
|B'_1(\ze)|&\le\f{\pi^2}2\sum_{n=N+1}^\infty\f{1-r_n^2}{(1-r_n)^2+(\ph_n-t)^2}\\
&\le\f{\pi^2}2\sum_{n=N+1}^\infty\f{1-r_n^2}{(1-r_n)^2+\ph_n^2}\\
&\le\f{\pi^2}2\sum_{n=N+1}^\infty\f{1-r_n^2}{|1-z_n|^2}<\f\eta4,
\end{aligned}
\end{equation}
where the last step relies on \eqref{eqn:choiceofn}. Eventually, we obtain 
\begin{equation}\label{eqn:bprimeeta}
|B'_1(\ze)|<\f\eta4,\qquad\ze\in\ga_\de
\end{equation}
(we have actually checked this for the bigger arc $\{\ze\in\T:\,-\f\pi2<\arg\,\ze<0\}$, not just for $\ga_\de$). 
Finally, we combine \eqref{eqn:fprimebelow} with \eqref{eqn:gongade} and \eqref{eqn:bprimeeta} to conclude that 
$$|f'(\ze)|\ge\f\eta4,\qquad\ze\in\ga_\de.$$ 
This yields the left-hand side inequality in \eqref{eqn:zhaba}, with the appropriate $C$, and completes the proof. 

\medskip

\section{Two examples}

The purpose of this section is to show that conditions \eqref{eqn:firstcond} and \eqref{eqn:secondcond} 
appearing in Theorem \ref{thm:newmainres}, via the definition of $\si_{\mathcal E}(f)$, are indispensable 
and close to being sharp. 
\par The two examples below follow the same pattern (and many more relevant examples can be furnished along 
these lines). Let $h:\D\to\Om$ be a conformal mapping of the disk onto a domain $\Om$ that is contained in 
the left half-plane 
$$\mathcal H:=\{w\in\C:\,\text{\rm Re}\,w<0\}$$ 
and contains, for some fixed number $c$, infinitely many 
points of the form $c+2\pi ik$ with $k\in\Z$. More precisely, we are assuming that there is a $c\in\mathcal H$ 
and an infinite subset $\La$ of $\Z$ such that 
$$\{c+2\pi ik:\,k\in\La\}\subset\Om.$$ 
(This is certainly the case if $\Om$ contains a vertical line or half-line.) Further, put $a:=e^c$ and note 
that $|a|<1$; then define 
\begin{equation}\label{eqn:violin}
g:=e^h\qquad\text{\rm and}\qquad f:=\f{g-a}{1-\bar ag}.
\end{equation}
Since $h(\D)=\Om\subset\mathcal H$, it follows that $g$ (and hence also $f$) is an $H^\infty$-function of norm 
at most $1$. Now, suppose that $h$ maps a certain arc $\Gamma\subset\T$ continuously onto an interval -- possibly 
infinite -- of the imaginary axis $i\R$. We have then $|g|=|f|=1$ on $\Gamma$ (whence $\|g\|_\infty=\|f\|_\infty=1$), 
and moreover, $g$ and $f$ will each have an angular derivative on $\Gamma$. In addition, $f$ vanishes at the 
points $z_k:=h^{-1}(c+2\pi ik)$ with $k\in\La$, because $g(z_k)=a$. Consequently, letting $B$ denote the Blaschke 
product with zeros $z_k$, $k\in\La$, we see that the inner part of $f$ is divisible by $B$. 
\par On the other hand, 
\begin{equation}\label{eqn:viola}
f'=\f{1-|a|^2}{(1-\bar ag)^2}gh'.
\end{equation}
The function $(1-\bar ag)^{-2}$ is outer (and even invertible in $H^\infty$); therefore, if $g$ and $h'$ also 
happen to be outer, the same will be true for $f'$. The situation then stands in sharp contrast to the conclusion 
of Theorem \ref{thm:newmainres}: indeed, $f$ has a nonconstant inner factor, while $f'$ has none. This means that 
the current function $f$ violates the hypotheses of the theorem. Specifically, if $\mathcal E$ is taken to be 
$\Gamma$ (so that $\widetilde{\mathcal E}=\T\setminus\Gamma$) and if the zero sequence $\{z_k\}$ is thick, then 
either \eqref{eqn:firstcond} or \eqref{eqn:secondcond} must break down. Thus, a glance at a concrete example of 
the above type might reveal whether the two sufficient conditions are reasonably close to being necessary. 

\medskip\noindent\textbf{Example 1.} Let 
$$\Om=\{w\in\C:\,-\pi<\text{\rm Re}\,w<0\}.$$ 
The function 
$$h(z)=i\log\left(\f{1+z}{1-z}\right)-\f\pi 2$$ 
(where the principal branch of the logarithm is used) maps $\D$ conformally onto $\Om$. We now fix a number $c$ 
with $-\pi<c<0$, then put $a=e^c$ and define the functions $g$ and $f$ by \eqref{eqn:violin}, with the current 
$h$ plugged in. Letting $\mathcal E$ stand for the arc $\{e^{it}:-\pi<t<0\}$, we have $h(\mathcal E)=i\R$, and 
so $f$ has an angular derivative on $\mathcal E$. Since $e^{-\pi}\le|g|\le1$ on $\D$, it follows that $g$ is outer. 
The function $h'(z)=2i(1-z^2)^{-1}$ being outer as well, we may invoke \eqref{eqn:viola} to deduce that $f'$ is outer. 
\par Now, the zeros $z_k$ of $f$, given by 
\begin{equation}\label{eqn:zedk}
z_k=h^{-1}(c+2\pi ik),\qquad k\in\Z,
\end{equation}
have the property that $\om_{z_k}(\widetilde{\mathcal E})$ takes the constant value $|c|/\pi$. 
(This is best seen by looking at the images $\ze_k$ of the $z_k$'s under 
the transformation 
\begin{equation}\label{eqn:cmdhp}
z\mapsto\f{1+z}{1-z},
\end{equation}
which maps $\D$ onto the right half-plane and $\mathcal E$ onto 
the half-line $i\R_-:=\{i\eta:\eta<0\}$. The points $\ze_k=(1+z_k)/(1-z_k)$ are then determined by the formula 
\begin{equation}\label{eqn:zetak}
\ze_k=\exp\left(2\pi k-ic-\f{i\pi}2\right),\qquad k\in\Z,
\end{equation}
whence 
\begin{equation}\label{eqn:anglezetak}
\arg\ze_k=-\f\pi 2-c=-\f\pi 2+|c|. 
\end{equation}
Recalling the well-known interpretation of the harmonic measure 
in terms of angles, one readily arrives at the required fact.) 
\par Finally, we observe that the sequence $\{z_k\}$ clusters at the points $\pm1$ and is thick. The latter claim 
can be verified with the help of a lemma by Sundberg and Wolff from \cite{SW}. (Precisely speaking, the version we 
need is obtained by combining Lemma 7.1 on p.\,578 of \cite{SW} with the concluding paragraph on p.\,580 after 
the lemma's proof. In fact, \cite{SW} treats a more general situation involving a Douglas algebra $B$, 
which we take to be $H^\infty+C$. See also \cite[p.\,4455]{DN} for the special case in question.) To state the 
thinness criterion given there, let $\{a_j\}$ be a sequence of distinct points in $\D$. Also, consider the arcs 
$$I_{N,j}:=\{\ze\in\T:\,|\ze-a_j|\le N(1-|a_j|)\}$$ 
with $N>1$, and write $\mathcal K(N,j)$ for the set of those indices $k$, $k\ne j$, which satisfy 
$a_k/|a_k|\in I_{N,j}$ and $1-|a_k|\le m(I_{N,j})$. This done, the Sundberg--Wolff result tells us that 
$\{a_j\}$ is thin if and only if, for every $N>1$, 
\begin{equation}\label{eqn:suwo}
\lim_{j\to\infty}(1-|a_j|)^{-1}\sum_{k\in\K(N,j)}(1-|a_k|)=0.
\end{equation}
Now, a computation shows that setting $a_j=z_j$, where the $z_j$'s are given by \eqref{eqn:zedk}, makes 
\eqref{eqn:suwo} false (provided that $N$ is large enough). Again, the easiest way to check this is to 
rephrase \eqref{eqn:suwo} for the right half-plane and then look at the images \eqref{eqn:zetak} of 
the $z_k$'s under the conformal mapping \eqref{eqn:cmdhp}. Thus, $\{z_k\}$ is indeed a thick sequence. 
\par In summary, while a suitably tangential convergence (in the sense of \eqref{eqn:firstcond} and 
\eqref{eqn:secondcond}) of the zero sequence $\{z_k\}$ to the endpoints $\pm1$ of $\mathcal E$ would 
imply that $f'$ has an inner factor, no kind of nontangential convergence would suffice. 
In fact, \eqref{eqn:anglezetak} shows that, letting $c$ be appropriately small in modulus, we can 
arrange it for the $\ze_k$'s to lie on a half-line that forms an arbitrarily small angle with the 
lower imaginary semiaxis (or, equivalently, for $z_k$ to lie on a circular arc that forms 
an arbitrarily small angle with $\mathcal E$ at $\pm1$). 

\medskip\noindent\textbf{Example 2.} Now let 
$$\Om=\{w\in\C:\,\text{\rm Re}\,w<0,\,\text{\rm Im}\,w<0\}.$$ 
This time, we take the conformal map $h:\D\to\Om$ to be 
$$h(z)=-e^{i\pi/4}\sqrt{\f{1+z}{1-z}},$$
where the square root is supposed to satisfy $\sqrt x>0$ for $x>0$. We then fix a number $c\in(-\infty,0)$ 
and define the functions $g$ and $f$ by \eqref{eqn:violin}, with $a=e^c$. 
\par This done, we claim that $g$ is an outer function. Indeed, since $g$ is zero-free and has radial limit $0$ 
only at the point $1$, it follows that $g$ has no inner factor, except possibly for the \lq\lq atomic" singular 
function 
$$S_\ga(z):=\exp\left(\ga\f{z+1}{z-1}\right)$$ 
with some $\ga>0$. However, if $g$ were divisible by $S_\ga$, then we would have 
$$|g(x)|\le\exp\left(\ga\f{x+1}{x-1}\right),\qquad0<x<1,$$
whereas $g$ actually has a milder decay rate as $x\to1^-$; in fact, 
$$|g(x)|=\exp\left(-\f1{\sqrt2}\cdot\sqrt\f{1+x}{1-x}\right),\qquad0<x<1.$$ 
Thus, $g$ is outer. So is the function 
$$h'(z)=-e^{i\pi/4}(1+z)^{-1/2}(1-z)^{-3/2},$$
and we eventually conclude, by virtue of \eqref{eqn:viola}, that $f'$ is outer also. 
\par The semicircle $\{e^{it}:0<t<\pi\}=:\mathcal E$ is mapped by $h$ onto the half-line $i\R_-$, 
so $f$ has an angular derivative on $\mathcal E$. (Note that the current $\mathcal E$ is different from 
its namesake in Example 1.) Finally, the zeros $z_k$ of $f$ are now given by 
$$z_k=h^{-1}(c-2\pi ik),\qquad k\in\N.$$ 
Equivalently, the points $\ze_k=(1+z_k)/(1-z_k)$ (i.e., the images of the $z_k$'s in the right half-plane 
under the transformation \eqref{eqn:cmdhp}) are determined by 
\begin{equation}\label{eqn:gandon}
\ze_k=-i(c-2\pi ik)^2,\qquad k\in\N.
\end{equation}
\par We have $z_k\to1$ and 
\begin{equation}\label{eqn:compar}
\om_{z_k}\left(\widetilde{\mathcal E}\right)\asymp\f1k,\qquad k\in\N.
\end{equation} 
To verify \eqref{eqn:compar}, one may first rewrite \eqref{eqn:gandon} in the form 
$$\ze_k=\xi_k+i\eta_k,$$
where 
$$\xi_k=4\pi|c|k\quad\text{\rm and}\quad\eta_k=4\pi^2k^2-c^2.$$ 
Now, the image of $\widetilde{\mathcal E}$ under the map \eqref{eqn:cmdhp} is $i\R_-$, and the angle 
at which this half-line is seen from $\ze_k$ is comparable to its tangent, $\xi_k/\eta_k$ 
(or equivalently, to $1/k$). Moving back to the disk, one arrives at \eqref{eqn:compar}. Furthermore, 
a calculation shows that 
$$1-|z_k|\asymp\f1{k^3},\qquad k\in\N.$$ 
Together with \eqref{eqn:compar}, this ensures that 
$$\om_{z_k}\left(\widetilde{\mathcal E}\right)\cdot\log\f1{1-|z_k|}\to0,$$
making \eqref{eqn:firstcond} true. By contrast, \eqref{eqn:secondcond} breaks down, the reason being 
that $|f'|$ becomes too small near the endpoint $1$ of $\widetilde{\mathcal E}$. 
\par The conclusion is that condition \eqref{eqn:firstcond} alone, or even its stronger version 
$$\om_{z_k}\left(\widetilde{\mathcal E}\right)=O\left((1-|z_k|)^{1/3}\right),$$
is not enough to guarantee the validity of our Gauss--Lucas type phenomenon (i.e., to ensure that 
$f'$ has a nontrivial inner factor whenever $f$ does). A different type of tangency condition, stated 
in terms of $|f'|$, should be added to make things work. 
\par Finally, we remark that the sequence $\{z_k\}$ in this last example was thick. This, again, can be 
verified by means of the Sundberg--Wolff criterion \eqref{eqn:suwo}, possibly transplanting everything 
to the right half-plane (for the sake of convenience) and working with the $\ze_k$'s instead.

\end{document}